\documentclass[a4paper,twoside]{article}
\usepackage[latin1]{inputenc}

\usepackage{amsmath,amsthm,amsfonts,latexsym,amscd,amssymb,enumerate}
\newtheorem{theorem}{Theorem}
\newtheorem{lemma}[theorem]{Lemma}
\newtheorem{corollary}[theorem]{Corollary}

\theoremstyle{definition}
\newtheorem{definition}[theorem]{Definition}

\theoremstyle{remark}
\newtheorem{remark}[theorem]{Remark}

\newcommand{\Boutet}{\mathfrak{A}}
\newcommand{\compacts}{\mathfrak{K}}
\newcommand{\ideal}{\mathfrak{I}}

\newcommand{\reals}{\mathbb{R}}
\newcommand{\complexs}{\mathbb{C}}
\newcommand{\naturals}{\mathbb{N}}
\newcommand{\integers}{\mathbb{Z}}

\DeclareMathOperator{\id}{id}
\newcommand{\boundary}[1]{\partial#1}
\newcommand{\boundedops}{\mathcal{B}}

\newcommand{\into}{\hookrightarrow}

\newcommand{\iso}{\cong}

\newcommand{\cA}{{\mathcal A}}
\newcommand{\fA}{{\mathfrak A}}
\newcommand{\fH}{{\mathfrak H}}
\newcommand{\fK}{{\mathfrak K}}
\newcommand{\fT}{{\mathfrak T}}
\newcommand{\intX}{X^\circ}
\newcommand{\lh}{{\mathfrak L}(\mathfrak H)}
\newcommand{\rp}{{\mathbb R}_+}

\DeclareMathOperator{\im}{im}      
\DeclareMathOperator{\End}{End}    
\DeclareMathOperator{\Hom}{Hom}    
\DeclareMathOperator{\Id}{Id}

\DeclareMathOperator{\ind}{ind}

\newcommand{\forget}[1]{}

\allowdisplaybreaks[2]

\begin{document}
\pagestyle{myheadings}
\markboth{S.~Melo, E.~Schrohe, T.~Schick}%
{K-theoretic Proof of Boutet de Monvel's Index Theorem}
\date{ }

\title{A K-theoretic Proof of Boutet de Monvel's\linebreak[4] 
Index Theorem for Boundary Value Problems}

\author{
Severino T. Melo
\and 
Thomas Schick
\and
Elmar Schrohe
}\maketitle

\begin{abstract}
We study the $C^*$-closure $\fA$ of the algebra of all operators 
of order and class zero in Boutet de Monvel's
calculus on a compact connected manifold $X$ with boundary 
$\partial X\not=\emptyset.$ 
We find short exact 
sequences in $K$-theory
\begin{equation*}
   0\to  K_i(C(X))\to K_i(\Boutet/\compacts) 
\stackrel p\to
   K_{1-i}(C_0(T^*X^\circ)) \to 0 , \quad i= 0,1,
 \end{equation*}
which split, so that $K_i(\Boutet/\compacts)\iso 
K_i(C(X))\oplus K_{1-i}(C_0(T^*X^\circ))$.
Using only simple $K$-theoretic arguments and the Atiyah-Singer Index Theorem, 
we show that the Fredholm index of an
elliptic element in $\cA$ is given by 
$$\ind A=\ind_t(p([A])),
$$
where $[A]$ is the class of $A$ in $K_1(\Boutet/\compacts)$
and $\ind_t$ is the topological index, 
a relation first established by Boutet de Monvel by different methods.  
\medskip

{\bf  Math.\ Subject Classification}:  58J32, 19K56, 46L80

\end{abstract}

\section*{Introduction}
\label{sec:introduction}

Boutet de Monvel's calculus 
is a pseudodifferential calculus on manifolds with boundary.
It comprises the classical differential boundary value problems
as well as the parametrices to elliptic elements, even their inverses whenever they exist \cite{MR53:11674}.
If the underlying manifold is compact, elliptic operators 
in the calculus define Fredholm operators between appropriate Hilbert spaces. 
Boutet de Monvel established an index theorem for that case:
He showed that there is a map which associates to each elliptic operator
an element in the $K$-theory  of the cotangent bundle over the interior of
the manifold and that the Fredholm index is the composition of 
that map with the ($\integers$-valued) topological index map. 

The crucial step is the construction
of the above map from elliptic operators to $K$-theory,
for which Boutet de Monvel uses elaborate considerations, combining
homotopy arguments within the algebra with classical
(vector-bundle) $K$-theory.
In this article we will show that this map 
can be obtained relying only on basic knowledge of the 
structure of the algebra and relatively simple constructions in 
$K$-theory for $C^*$-algebras  (which were not yet available in 1971). 
Boutet de Monvel's map is neither obvious nor trivial.
The point is that we are able to represent it as a composition of 
various standard (yet non-trivial) maps in $K$-theory.

To be more specific, 
let $X$ be a compact $n$-dimensional manifold with boundary $\partial X$, 
embedded in a closed manifold $\widetilde X$ of the same dimension
($\tilde X$ could e.g.~be the double of $X$). 
By $\intX$ we denote the interior of $X$.
We assume that $X$ is connected and $\partial X$ is nonempty.

Given a pseudodifferential operator $P$ on $\widetilde X$, we define
the truncated pseu\-do\-dif\-ferential operator 
$P_+:C^\infty(X)\to C^\infty(\intX)$ 
as the composition
$r^+Pe^+$, where $e^+$ is extension by zero from $X$ to $\widetilde X$ and 
$r^+$ is the restriction to $\intX$. In general, the functions in the 
range of $P_+$ will not be smooth up to the boundary. 
One therefore assumes that $P$ satisfies the {\em transmission
condition}, a condition on the symbol of $P$ 
which ensures that both $P_+$ and $(P^*)_+$, the truncated 
operator formed from the formal adjoint of $P$, map smooth 
functions on $X$ to smooth functions on $X$. 

An operator in Boutet de Monvel's calculus is a matrix
\begin{eqnarray}
\label{1}
A=
\begin{pmatrix}P_++G&K\\T&S
\end{pmatrix}: 
\begin{array}{ccc}
C^\infty(X,E_1)& & C^\infty(X,E_2)\\
\oplus&\to&\oplus\\
C^\infty(\partial X,F_1)&&C^\infty(\partial X,F_2)
\end{array} 
\end{eqnarray} 
acting on sections of vector bundles $E_1$, $E_2$ over $X$  and $F_1$, $F_2$ over $\partial X$.
Here, $P$ is a pseudodifferential operator satisfying the transmission condition; 
$G$ is a singular Green operator, 
$T$ is a trace operator, 
$K$ is a potential (or Poisson) operator, 
and $S$ is a pseudodifferential operator on $\partial X$. 
All these operators are assumed to be classical; i.e.\ their symbols have polyhomogeneous expansions 
in the respective classes.
For details, the reader is referred to the monographs
by Rempel and Schulze \cite{RS} or Grubb \cite{G} as well as to the short introduction \cite{S3}.
We will need the following facts:

The operators $G$, $K$, and $T$\
are regularizing in the interior of $X$.
In a collar neighborhood of the boundary, they can be viewed as 
operator-valued pseudodifferential operators along the boundary. 
In particular, they have an order assigned to them.
The singular Green and the trace operators 
also have a {\em class}\ (or {\em type}) $d\in\naturals_0$, 
related to the order of the derivatives appearing in the 
boundary condition. 

The composition of two operators  of the form \eqref{1}
is defined whenever the vector bundles serving as the range of the first operator 
form  the domain of the second.
The composition $AA'$ of an operator $A'$ of order $m'$ and class 
$d'$ with an operator $A$ of order $m$ and
class $d$ results in an operator of order $m+m'$ and  class $\le \max(m'+d, d')$. 
In particular, the composition of two operators of order and 
class zero is again of order and class zero.

For $E_1=E_2=E$ and $F_1=F_2=F$, the operators of order and class 
zero thus form an algebra $\cA$. 
Moreover, they extend to bounded operators on the 
Hilbert space $\fH=L^2(X,E)\oplus L^{2}(\partial X,F)$.
In fact, $\cA$ is a $*$-subalgebra of $\lh$,
closed under holomorphic functional calculus, cf.\ \cite{S2}.
We use here the definition of order and class in \cite{RS} and \cite{S3} 
which differs slightly from that in \cite{G}.
It allows us to use the $L^2$-space over the boundary instead of $H^{-1/2}$
and gives us better homogeneity properties of the boundary symbols.
In view of the fact that both the kernel and the cokernel of an elliptic 
operator in $\cA$ consist of smooth functions, the choice is irrelevant 
for index theory.

Standard reductions - recalled in Section \ref{sec:reduct-index-probl} -
allow to reduce any index problem to the case where the
operator is an element of the algebra $\cA$, 
so that we can apply operator-algebraic methods.
This is a central point of the paper.

Let us have a closer look at the structure of $\cA$.
In a generalization of the classical Lopatinskij-Shapiro condition, 
the ellipticity of an element $A\in \cA$ is governed by
the invertibility of two symbols, 
namely the pseudodifferential principal symbol, $\sigma(A)$, 
and the principal boundary symbol, $\gamma(A)$, which take 
values in certain $C^*$-algebras over $S^*X$ and $S^*\partial X$, 
respectively, cf.\ Section \ref{sec:symbol-boutet-de}.

The maps 
$\sigma$ and $\gamma$ are $*$-homomorphisms on $\cA$.
Extending the classical results by Gohberg and Seeley, 
Rempel and Schulze (\cite{RS}, 2.3.4.4, Theorem~1,  
based on work by Grubb and Geymonat \cite{GG})
showed that 
\begin{equation}
\label{rs}
\inf_{C\in\fK}||A+C||=\max\{||\sigma(A)||,||\gamma(A)||\}, \ \text{for all}\ A\in\cA,
\end{equation}
where $\fK$\ denotes the ideal of the compact operators  on $\fH$, 
and the norms on the right hand side are the supremum norms on $S^* X$
and $S^*\partial X$, respectively. 

We shall now denote by $\fA$ the closure (equivalently, the  $C^*$-closure)  
of $\cA$ in the topology of $\lh$. 
It follows from \eqref{rs} that  $\sigma$ and $\gamma$ 
have continuous extensions to $\fA$; we denote them  by the same letters.
An element of $\fA$ is compact,  if and only if  both symbols are zero.
Moreover, it is a Fredholm operator if and only if it is elliptic, i.e.,
both symbols are (bundle) isomorphisms. 

Boutet de Monvel showed that, 
given  an elliptic element $A$ 
in $\cA$, one can 
find a homotopy through elliptic elements in $\cA$,
connecting $A\oplus \Id$ to an operator of the 
form $\begin{pmatrix}P'_+&0\\0&Q'\end{pmatrix}$,
where $P'$ is a pseudodifferential operator whose principal symbol 
is an isomorphism of $E$ (independent of $\xi$) in a neighborhood
of the boundary.  
Through the usual difference bundle construction,
the principal symbol of $P'$ therefore defines an element 
$[P']$ of $K(T^*\intX)$. Boutet de Monvel then associated to 
$A$ the class $[P']+ {\rm Thom}([Q'])$, 
where $[Q']$ is the class in $K(T^*\partial X)$ defined by the 
principal symbol of $Q'$ and 
${\rm Thom}:K(T^*\partial X)\to K(T^*\intX)$
is the Thom map, also called Gysin homomorphism, or Umkehrmap. 
Moreover, he showed that the composition of that map with the topological 
index map $K(T^*\intX)\to \integers$ gives the index of $A$.

Fedosov \cite{MR1401125} then proved that this implies the formula 
\begin{equation}\label{Fedosov}
\ind A = \int_{T^*X} {\rm ch}(\sigma(A)){\mathcal T}(X) + \int_{T^*\partial X}{\rm ch}'(\gamma(A)){\mathcal T}(X).
\end{equation} 
Here ${\mathcal T}(X)$ is the differential form given by the Todd
class of the tangent bundle (or its restriction to the boundary), 
${\rm ch}\,\sigma(A)$ is the Chern character associated with the $K$-class 
induced by the pseudodifferential principal symbol 
$\sigma(A)$, and ${\rm ch'}\gamma(A)$ is a differential form constructed from the 
$K$-class of the boundary symbol; it is given by a formula analogous to that for the 
usual Chern character.

In order to establish this formula, Fedosov referred to 
Boutet de Monvel's work and showed two facts:
(i) 
the formula is invariant under homotopies 
within the class of elliptic boundary value problems in the calculus and 
(ii) 
whenever the principal symbol of $P$ is an isomorphism in a neighborhood of the boundary in $X$,
independent of the covariable $\xi$, and the boundary symbol is an isomorphism over the full ball
bundle $B^*(\partial X)$, then the above formula reduces to the classical formula of Atiyah and Singer. 

Our approach to the index theorem is based on a careful analysis of the boundary symbol map. 
Indeed, since the kernel of $\gamma$ contains the compact ope\-rators, we have a natural short exact sequence
\begin{equation}
\label{sesgamma}0\to \ker \gamma/\fK \to \fA/\fK \xrightarrow{\gamma} \im\gamma \iso\fA/\ker \gamma\to 0.\end{equation}
What is mainly needed for the understanding of our proof of the index theorem is 
the structure of $\ker\gamma$ and $\im \gamma.$
Both were determined in \cite[Section 3]{MR1998610} for the case of trivial one-dimensional bundles $E$ and $F$.
We shall review these computations in Section 1 for the case of general bundles.

In Section 2, we recall a  basic construction from $K$-theory, namely how a commutative diagram of 
short exact sequences of $C^*$-algebras 
\begin{equation*}
    \begin{CD}
      0 @>>> A @>>> B @>>> C @>>> 0\\
      &&@AA{f}A @AA{g}A @AA{h}A\\
      0 @>>> A' @>>> B' @>>> C' @>>> 0
    \end{CD}
  \end{equation*}
yields a commutative grid of $C^*$-algebras involving $A', B', C'$, the mapping cones,
and the suspensions of $A,B,$ and $C$.
We then apply this in Section 3 to the above sequence \eqref{sesgamma}, linked to the sequence 
$$
0\to C_0(\intX)\to C(X)\to C(\partial X)\to 0. 
$$ 
As a first result, we obtain 
\begin{theorem}\label{t1}
 We have natural short exact sequences
  \begin{equation}\label{eq:main}
   0\to  K_i(C(X))\longrightarrow K_i(\Boutet/\compacts) \stackrel{p}{\longrightarrow}
   K_{1-i}(C_0(T^*X^\circ)) \to 0 ,
 \end{equation}
$i=0,1$, which split; i.e.~we have  $($not necessarily natural$)$ isomorphisms
$$
  K_i(\Boutet/\compacts) \iso K_i(C(X)) \oplus K_{1-i}(C_0(T^*X^\circ)).
$$
\end{theorem}
This strengthens the results of Melo, Nest, and Schrohe \cite[Corollaries 12 and 19]{MR1998610}, 
where the corresponding statements 
were derived for $i=0$ under an additional hypothesis and for $i=1$ 
using Boutet de Monvel's index theorem. 

In Section 4 we then prove the index theorem:
\begin{theorem}\label{t2}
The index of a Fredholm operator $A$ in $\Boutet$
is given by 
\begin{equation}\label{bdmif}
\ind A=\ind_t(p([A])),
\end{equation}
where $[A]$ is the $K_1$-class of  
$A$ in $\Boutet/\compacts$, $p$ is the map in 
\eqref{eq:main} and 
$$\ind_t: 
 K_{0}(C_0(T^*X^\circ))\iso K^0(T^*X^\circ))\to K^0(pt)\iso \mathbb Z$$
is the topological index. Moreover, $p[A]={\tt ind}(A)$, 
where ${\tt  ind}$ denotes the map defined by Boutet de Monvel in {\rm \cite[Theorem 5.21]{MR53:11674}}.
\end{theorem}

As in the classical case, Fedosov's arguments yield the cohomological form \eqref{Fedosov} of the index theorem. 

\section{Elliptic Operators and Symbols in Boutet de Monvel's Calculus}
\label{sec:preliminaries}

\subsection{A Normal Form for the Index Problem}
\label{sec:reduct-index-probl}
Suppose we are given a Fredholm operator $A$ in Boutet de Monvel's calculus of 
oder $m$ and class $d$ acting on sections of vector bundles over $X$  as in \eqref{1}.

There exist order and class reducing 
invertible operators in the
calculus. By invertibility, composition with those does 
not change the index. 
Therefore, we can always achieve that order and class are zero.
 
Moreover, we can always assume that $X$ is connected. If it is not,
then the spaces $C^\infty(X,E_1)$, \ldots, decompose as direct sums
corresponding to the decomposition in connected components, 
and the operator $A$ becomes a matrix $M$ with respect to this
decomposition. By definition of the Boutet de Monvel calculus, the
off-diagonal entries  
are smoothing operators, thus compact. 
The index is therefore unchanged if we replace
$M$ by the diagonal matrix equal to the diagonal of
$M$. But then it is clear that the index is additive, and the
individual entries correspond to the connected components of $X$.

We can also assume that $E_1=E_2=E$ and $F_1=F_2=F$. 
Indeed, if an elliptic Boutet de
Monvel operator between $(E_1,F_1)$ and $(E_2,F_2)$ exists, we have in
particular an interior elliptic pseudodifferential operator $Q$ between
$E_1$ and $E_2$. By definition, its symbol defines an isomorphism
between $\pi^*E_1$ and $\pi^*E_2$, where $\pi\colon S^*X\to X$ is the
projection of the unit cotangent sphere bundle. Since $X$ has a
non-empty boundary, $\pi\colon S^*X\to X$ has a section, 
see e.g.\ \cite[Proposition 9]{MR1998610} 
or use that the Euler class of the cotangent bundle is
trivial. 
Restricting
the symbol isomorphism to this section, we get a bundle isomorphism
between $E_1$ and $E_2$. 
 
Next we choose a complement $\widetilde E$ to $E$ such that 
$E\oplus \widetilde E=\underline N$, the $N$-dimensional trivial
bundle. We can take $N$ so large that, over $\partial X$, the
bundles $F_1$ and $F_2$ are also embedded in the restriction of
$\underline N$ 
and consider the operator 
$$\widetilde A=\begin{pmatrix}\id&0\\0&A\end{pmatrix}:
\begin{array}{ccc}
C^\infty(X,\widetilde E)& & C^\infty(X,\widetilde E)\\
\oplus&&\oplus\\
C^\infty(X,E)& \to& C^\infty(X,E)\\
\oplus&&\oplus\\
C^\infty(\partial X,F_1)&&C^\infty(\partial X,F_2)
\end{array},  $$
which has the same index as $A$.
 
We then recall that there are elliptic operators $R_j$, $j=1,2$, of order and class
zero in Boutet de Monvel's calculus such that 
$$R_j=
\begin{pmatrix}
\Gamma_{F_j} \\ p_{F_j}\circ \gamma_0
  \end{pmatrix}:C^\infty(X,\underline N)\to 
\begin{array}{c}C^\infty(X,\underline N)\\\oplus\\  
C^\infty(\partial X,F_j)\end{array}$$ is a Fredholm operator of index zero,
\cite[Theorem (5.12)]{MR53:11674}).
Here
$\gamma_0$ is restriction to the boundary; $p_{F_j}$ is projection onto the
subbundle $F_j$.
Composing $\widetilde A$ from the left with a parametrix to $R_2$ and from the 
right with $R_1$ we obtain an operator with the same index as $A$ which is an endomorphism
of $C^\infty(X,\underline N)$.  

\subsection{Symbols}
\label{sec:symbol-boutet-de}
We consider an operator $A$ as in  \eqref{1}, with 
$E_1=E_2=E$, $F_1=F_2=F$. 
The pseudodifferential principal symbol $\sigma(A)$ of $A$ is defined to be the 
principal symbol of $P$, restricted to $S^*X$. 
The principal boundary symbol of $A$ is a smooth section from $S^*\partial X$ into
the endomorphisms of 
\begin{equation}\label{end}
\pi_\partial^*\left(L^2(\rp)\otimes E|_{\partial X}\right)\oplus \pi_\partial^*F\iso 
\left( L^2(\rp)\otimes \pi_\partial^*E|_{\partial X}\right)\oplus\pi_\partial^* F,
\end{equation}
where $\pi_\partial:S^*\partial X\to \partial X$ is the canonical projection.
It is best described for a trivial one-dimensional bundle and 
in local coordinates $(x',x_n,\xi',\xi_n)$ for 
$T^*X$ in a neighborhood of the boundary.
Here, $G$ acts like a pseudodifferential operator along
the boundary, with an operator-valued symbol taking values 
in regularizing operators in the normal direction.
One way to write this operator-valued symbol is via a so-called 
symbol kernel
$\tilde g = \tilde g(x',\xi',x_n,y_n)$. 
For fixed $(x',\xi')$, this is a rapidly decreasing function 
in $x_n$ and $y_n$ which acts as an integral operator on $L^2(\rp)$. 
It satisfies 
special estimates, combining the usual pseudodifferential
estimates in $x'$ and $\xi'$ with those for rapidly decreasing functions in
$x_n$ and $y_n$. 
The singular Green symbol $g$ of $G$ is defined from the symbol kernel via 
Fourier and inverse Fourier transform:
$$
g(x',\xi',\xi_n,\eta_n) = F_{x_n\to \xi_n}\overline F_{y_n\to \eta_n}\tilde
g(x',\xi',x_n,y_n).
$$
It has an expansion into homogeneous terms; the leading one we call $g_0$. Inverting
the operation above, we associate with $g_0$\ a symbol-kernel 
$\tilde g_0(x',\xi',x_n,y_n)$\ which is rapidly decreasing in $x_n$ and $y_n$
for fixed $(x',\xi')$. 
We denote by $g_0(x',\xi',D_n)$ the (compact) operator
induced on $L^2(\rp)$ 
by this kernel. 
Similarly,  $K$ and $T$ have symbol-kernels $\tilde k(x',\xi',x_n)$ and 
$\tilde t(x',\xi',y_n)$; these are  rapidly decreasing 
functions for fixed $(x',\xi')$. 
The symbols $k$ and $t$ are defined as their Fourier and inverse 
Fourier transforms. They have asymptotic expansions with leading terms
$k_0$ and $t_0$. Via the symbol-kernels $\tilde k_0$ and $\tilde t_0$
one defines $k_0(x',\xi', D_n): \complexs\to L^2(\rp)$  
as multiplication by $\tilde k_0(x',\xi',\cdot)$, while 
$t_0(x',\xi',D_n) :L^2(\rp)\to \complexs$ 
is the
operator $\varphi \mapsto \int \tilde t_0(x',\xi',\cdot)\varphi$.   

We denote by $p_0$ and $s_0$ the principal symbols of $P$ and $S$, respectively.
The boundary symbol $\gamma(A)$ of $A$ in $(x',\xi')$ is then defined by 
$$
\gamma(A)(x',\xi')=
\left(\begin{array}{cc}
p_0(x^\prime,0,\xi^\prime,D_n)_{_{+}}+g_0(x',\xi',D_n)&k_0(x',\xi',D_n)\\
t_0(x',\xi',D_n)&s_0(x',\xi')
\end{array}\right).
$$  
Two things are important to note: 
\begin{enumerate}\item 
Except for $p_0(x',0,\xi',D_n)$ all entries in $\gamma(A)(x',\xi')$ are compact.
\item 
The boundary symbol is `twisted' homogeneous of degree zero in the sense that 
$$\begin{pmatrix}\kappa_{\lambda^{-1}}&0\\0&\id\end{pmatrix}
\gamma(A)(x',\lambda\xi')
\begin{pmatrix}\kappa_{\lambda}&0\\0&\id\end{pmatrix}=\gamma(A)(x',\xi'),\quad\lambda>0,
$$
with the $L^2(\reals_+)$-unitary $\kappa_\lambda$ given by $\kappa_\lambda f(t)=\sqrt\lambda f(\lambda t)$.
\end{enumerate}

\subsection{Kernel and Range of the Boundary Symbol Map}
  \label{sec:kern-range-bound}

The algebra $\fA$ contains the ideal $\ideal$ given by the $C^*$-closure 
of all elements of  the form 
\begin{eqnarray*}
 \begin{pmatrix}\varphi P \psi +G&K\\T&S
\end{pmatrix}
\end{eqnarray*} 
with $\varphi$, $\psi$ in $C^\infty_c(X^\circ)$ and $G, K, T, S$
of negative order and class zero.
Clearly, $\ideal$ is contained in the kernel of $\gamma$. More is true:
\begin{theorem}\label{kergamma} The kernel of the boundary symbol map $\gamma$ is equal to $\mathfrak I$. 
The quotient $\mathfrak I/\compacts$ is isomorphic to $C_0(S^*X^\circ, \End \pi^*E)$ 
with isomorphism induced by the principal symbol.
Here $\pi^*E$ is the pull-back of  $E$ under the projection $\pi: S^*X^\circ\to X^\circ$.  
\end{theorem}
\begin{proof} 
This is immediate  from the considerations
for the case of trivial bundles \cite[Theorem 1]{MR1998610}.
\end{proof}

In order to make the computation of the range of $\gamma$ more transparent, let us first
consider the localized situation with $E$ and $F$ trivial one-dimensional.
We write $\gamma$ as a $2\times2$-matrix with entries $\gamma_{ij}$, $i,j=1,2$.

Let $p$ be a classical pseudodifferential symbol of order zero on $\reals^n$. 
For fixed $(x',\xi')$, $p_0(x',0,\xi',\xi_n)$ is a symbol of order zero on $\reals.$
The transmission property assures that the values of $p$ in $\xi_n=+\infty$ and
$\xi_n=-\infty$ coincide. 
The operator 
$$p_0(x',0,\xi',D_n)_+=r^+{\rm op}\,p_0(x',0,\xi',\xi_n)e^+:L^2(\reals_+)\longrightarrow L^2(\reals_+)$$
in the upper left corner 
$\gamma_{11}$ then is a Toeplitz type operator. 
In fact, it is unitarily equivalent to the usual Toeplitz operator $T_f$ with 
symbol $f(z)=p_0(x',0,\xi',i(z-1)/(z+1))$.    
Thus the image of the upper left corner under the boundary symbol map is a subalgebra of 
$C(S^*\partial X,\fT)$, where $\fT$ denotes the Toeplitz algebra.

All other entries in the matrix for $\gamma(A)(x',\xi')$ are compact, so that 
the boundary symbol is, for fixed $(x',\xi')$, a so-called Wiener-Hopf operator on 
$L^2(\reals_+)\oplus \complexs$. 
One might conjecture that the range of the boundary symbol map consisted of 
all sections with values in Wiener-Hopf operators. 
It came as a surprise (and turned out to be a crucial fact) 
in \cite{MR1998610} that this is not the case. 
It is the range of the upper left corner 
$\gamma_{11}$  which is slightly smaller than expected: 
Let us denote by $\fT_0$ the subalgebra of those Toeplitz operators whose symbol 
vanishes in $z=-1$ (corresponding to $\xi_n=\pm\infty)$.
The range of $\gamma_{11}$  
contains as an ideal all sections of $S^*\partial X$ with values in $\fT_0$, but
the only sections of the form $g(x',\xi')\otimes Id_{L^2(\reals_+)}$ it contains are
those where $g$ is independent of $\xi'$, thus a function on $\partial X$, not
$S^*\partial X$. We therefore get a split short exact sequence of
$C^*$-algebras
\begin{equation*}
  0\to C(S^*\boundary X, \fT_0) \to 
\im \gamma_{11} \to  C(\boundary X)\to 0,
\end{equation*}

Let us now go over to the case of general bundles, considering 
the entries in the matrix for $\gamma(A)(x',\xi')$ separately, writing $\tilde E$ and
$\tilde F$ instead of $\pi_\partial^*(E|_{\partial X})$  and $\pi_\partial^*F$:

\begin{enumerate}

\item 
The boundary symbol
$t_0(\cdot,\cdot,D_n)$ is a continuous 
section from $S^*\partial X$ to $\Hom(L^2(\reals_+)\otimes \tilde E,\tilde F)$ given by integration 
against the symbol kernel of $t_0$, hence a section of $\boundedops(L^2(\reals_+),\complexs)\otimes
\Hom(\tilde E,\tilde F)$. The construction in \cite[Lemma 4]{MR1998610}
shows that all elements in that space are obtained that way. 
\item Similarly, the range of the closure of the boundary symbol map for the Poisson
operators of order zero consists of all continuous sections from $S^*\partial X$ 
to $\boundedops(\complexs,L^2(\reals_+))\otimes \Hom(\tilde F,\tilde E)$.  
\item The boundary symbols of zero order pseudodifferential operators 
along the boundary are their principal symbols,
thus certain elements of $C(\partial X, \End \tilde F)$, and in fact, all elements in this
space are obtained as such symbols.
\item 
The continuous sections from $S^*\partial X$ into 
$\pi_\partial^*({\mathfrak T}_0\otimes \End(E|_{\partial X}))\iso 
{\mathfrak T}_0\otimes \End \tilde E $ are contained as an ideal in the range of 
$\gamma_{11}$ (the upper left corner of $\gamma$) 
by a bundle valued analog of \cite[Lemma 5]{MR1998610}, while, 
as in  \cite[Lemma 6]{MR1998610}, 
${\im} \gamma_{11}\cap C(S^*\partial X,\End \tilde E)
=C(\partial X, \End E|_{\partial X})$. 
Here, we consider the elements of $C(\partial X, \End E|_{\partial X})$ and
 $C(S^*\partial X, \End\tilde E)$ as elements of 
$C(S^*\partial X, \End(L^2(\reals_+)\otimes \tilde
E))$  by acting as the identity on $L^2(\reals_+)$.

We conclude that we get a split short exact sequence of $C^*$-algebras 
\begin{equation*}
  0\to C(S^*\partial X,{\mathfrak T}_0\otimes \End \tilde E) \to
  \im\gamma_{11} 
 \to C(\partial X, \End E|_{\partial X})\to 0.
\end{equation*}
\end{enumerate} 

We can now define the subbundle $\widetilde {\mathfrak W}_0$ 
of  endomorphisms of \eqref{end}, 
consisting of all $2\times 2$ matrices  
$w=(w_{ij})_{i,j=1,2}$
\forget{$\begin{pmatrix}w_{11}&w_{12}\\w_{21}&w_{22}\end{pmatrix} $}
 with $w_{11}\in {\mathfrak T}_0\otimes \End \tilde E$, $w_{12}\in \boundedops(\complexs,L^2(\reals_+))\otimes
\Hom(\tilde F,\tilde E)$, $w_{21}\in \boundedops(L^2(\reals_+),\complexs)\otimes \Hom(\tilde E,\tilde
F)$, and $w_{22}\in \End \tilde F$. We obtain:

\begin{theorem}\label{imgamma}
 The image of $\gamma$ fits into the following split exact sequence of
 $C^*$-algebras
 \begin{equation}\label{eq:split_ex_seq}
   0\to C(S^*\partial X, \widetilde {\mathfrak W}_0)\to 
\im\gamma\to
   C(\partial X, \End (E|_{\partial X}))\to 0.
 \end{equation}
\end{theorem}

We note that $\widetilde{\mathfrak W}_0$  is the bundle valued 
analog of the algebra  ${\mathfrak W}_0$ in \cite{MR1998610}; 
for $E=F=\complexs$ both coincide.
Strong Morita equivalence (as discussed in \cite[Section 1.5]{MR1998610}) 
together with the fact that $\mathfrak T_0$ has vanishing $K$-theory therefore implies 
(cf.~\cite[Lemma 7]{MR1998610})  

\begin{lemma}\label{KW0} $K_i(C(S^*\partial X, \widetilde{\mathfrak W}_0))=0$, $i=0,1.$
\end{lemma}

The split in \eqref{eq:split_ex_seq} is implemented by the $C^*$-algebra homomorphism 
$$b\colon C(\partial X,\End (E|_{\partial X}))\to {\im}\gamma,\quad 
g\mapsto  \gamma\left(\begin{pmatrix}f&0\\0&0\end{pmatrix}\right),$$ 
where $f$ is any continuous section in  $C(X,\End E)$ with $f|_{\partial X}=g$. 
We then conclude as in \cite[Corollary 8]{MR1998610}:

\begin{corollary}\label{biso} The induced homomorphism is an isomorphism
$$b_*\colon K_i(C(\partial X,\End(E|_{\partial X})))\to K_i({\im} \gamma)
=K_i(\Boutet/\ideal).$$
\end{corollary}

\section{K-theory Preliminaries}\label{sprl}

\begin{definition}
  Let $A$ be a $C^*$-algebra. The cone over $A$ is the $C^*$-algebra
  $CA:=\{f\colon [0,1]\to A\mid f(1)=0\}$. \end{definition}

Since $CA$ is a contractible $C^*$-algebra, its $K$-theory vanishes. The suspension of $A$ is given by     
$SA:= \{f\in CA\mid f(0)=0\}$.

\begin{definition}
  If $f\colon B\to A$ is a $C^*$-algebra homomorphism, the mapping
  cone $Cf$ is defined to be $Cf:=\{ (b,\phi)\in B\oplus CA;\;
  f(b)=\phi(0)\}$.
\end{definition}

 Projection onto $B$ defines a short exact sequence
  \begin{equation}\label{eq:cone_seq}
    0\longrightarrow SA{\mathop{\longrightarrow}\limits^{\hat\iota}} Cf{\mathop{\longrightarrow}\limits^{q}}B\longrightarrow 0.
  \end{equation}
The assignment of this exact sequence to each C$^*$-algebra homomorphism $f$ defines a functor between 
the corresponding categories (whose morphisms consist of commutative diagrams of homomorphisms or
of exact sequences, respectively). This functor is exact; i.e. we have:

\begin{lemma}\label{exf}
  Assume that the following is a commutative diagram of short exact
  sequences of $C^*$-algebras:
  \begin{equation*}
    \begin{CD}
      0 @>>> A @>>> B @>>> C @>>> 0\\
      &&@AA{f}A @AA{g}A @AA{h}A\\
      0 @>>> A' @>>> B' @>>> C' @>>> 0
    \end{CD}
  \end{equation*}
  Then we get an induced commutative grid of short exact sequences of
  $C^*$-algebras
  \begin{equation*}
    \begin{CD}
      && 0 && 0 && 0\\
      && @AAA @AAA @AAA\\
      0 @>>> A' @>>> B' @>>> C' @>>> 0\\
      && @AAA @AAA @AAA\\
      0 @>>> Cf @>>> Cg @>>> Ch @>>> 0\\
      && @AAA @AAA @AAA\\
      0 @>>> SA @>>> SB @>>> SC @>>> 0\\
      && @AAA @AAA @AAA\\
      && 0 && 0 && 0
    \end{CD}
  \end{equation*}
\end{lemma}

Lemma \ref{exf} can be proven by a diagram chase, using that the maps $CB\to CC$ and $SB\to SC$ are surjective (the exactness of $S$ is proven in \cite[Proposition 10.1.2]{R}). 

\begin{lemma}\label{l2}
The exact sequence \eqref{eq:cone_seq} induces six-term cyclic exact sequences in $K$-theory, whose connecting mappings 
$K_{i}(B)\to K_{1-i}(SA)$\ become, under the canonical isomorphisms $K_{1-i}(SA)\xrightarrow{\iso} K_i(A)$, the mappings induced by $f$. 
\end{lemma}

More precisely, this lemma states that the two diagrams
$$
\def\mapup#1{\Big\uparrow\rlap{$\vcenter{\hbox{$\scriptstyle#1$}}$}}
\begin{array}{ccc}
K_1(B)&{\mathop{\longrightarrow}\limits^{\delta_0}}&K_0(SA)\\
\mapup{=}&&\mapup{\Theta_A}\\
K_1(B)&{\mathop{\longrightarrow}\limits^{f_*}}&K_1(A)\end{array}
\ \text{and}\ 
\begin{array}{ccc}
K_0(B)&{\mathop{\longrightarrow}\limits^{\delta_1}}&K_1(SA)\\
\mapup{=}&&\mapup{\beta_A}\\
K_0(B)&{\mathop{\longrightarrow}\limits^{f_*}}&K_0(A)
\end{array}
$$
commute, where $\delta_0$ and $\delta_1$ denote, respectively, the index and the exponential mappings \cite[9.1.3 and 12.1.1]{R}
induced by \eqref{eq:cone_seq}, $\Theta_A$ is the isomorphism defined in \cite[10.1.3]{R},
and $\beta_A$ is the Bott isomorphism \cite[11.1.1]{R}. This follows
by applying  naturality of the long exact sequences
in $K$-theory to the diagram
\begin{equation*}
  \begin{CD}
    0 @>>> SA @>>> Cf @>>> B @>>> 0\\
 &&    @VVV @VVV @VV{f}V \\
    0@>>> SA @>>> CA @>>> A @>>> 0.
  \end{CD}
\end{equation*}

\begin{lemma}\label{l3} If $f:B\to A$ is a surjective C$^*$-homomorphism, then the map
$j:\ker f\ni x\mapsto (x,0)\in Cf$
induces a $K$-theory isomorphism, which fits into the commutative diagram
$$
\def\mapup#1{\Big\uparrow\rlap{$\vcenter{\hbox{$\scriptstyle#1$}}$}}
\begin{array}{ccccccccc}
\longrightarrow&K_{i+1}(B)&\longrightarrow&K_i(SA)&{\mathop{\longrightarrow}\limits^{\hat\iota_*}}&K_i(Cf)&
{\mathop{\longrightarrow}\limits^{q_*}}&K_i(B)&\longrightarrow\\
                &\mapup{=}&     &\mapup{\delta_{i+1}}&          &\mapup{j_*}&             &\mapup{=}&\\
\longrightarrow&K_{i+1}(B)&{\mathop{\longrightarrow}\limits^{f_*}}&K_{i+1}(A)&\longrightarrow&K_i(\ker f)&\longrightarrow&K_i(B)&\longrightarrow
\end{array}
$$
where the upper row is the cyclic exact sequence induced by \eqref{eq:cone_seq}, and the lower one is that induced by
$$
0\longrightarrow \ker f\longrightarrow B{\mathop{\longrightarrow}\limits^{f}} A\longrightarrow 0.
$$
\end{lemma}
\begin{proof}
  This result is certainly well known. For the sake of completeness,
  and since we did not find a convenient reference, we sketch the
  arguments. 
  One proves that $j_*$ is a $K$-theory isomorphism using that
  $K_*(CA)=0$ in the cyclic exact sequence associated to the short
  exact sequence
\[
0\longrightarrow \ker f\ {\mathop{\longrightarrow}\limits^{j}}\ Cf
\longrightarrow CA\longrightarrow 0
\] 
induced by projection of $Cf$ onto the second coordinate.  The boundary
map
$\delta_{i+1}$ is well known to be an isomorphism. It remains to
establish commutativity. The commutativity of the
left rectangle is part of Lemma \ref{l2}. The right rectangle commutes
by naturality of the $K$-theory functor, since $q\circ j= (i
\colon
\ker f\into B)$. The argument for the middle square is a little more
involved. Observe that $j_*^{-1} \hat\iota_* \delta_{i+1}$ defines another
homomorphism $K_{i+1}(A)\to K_i(\ker f)$ which is, by the naturality of
all constructions, natural and makes the $K$-theory sequence of the short
exact sequence of $C^*$-algebras exact.

However, homomorphisms with this properties are defined uniquely (up to
a universal sign) \cite[Exercise 9.F]{WO}, therefore the diagram is commutative up to this
sign. The special exact sequence $0\to SA\to CA\to A\to 0$ shows that
this sign is $+1$.
\end{proof}

\section{K-theory of Boutet de Monvel's Algebra}\label{sprf}

In order to keep the notation simple, we shall write $C_0(X^\circ)$, $C(X), $ and $C(\partial X)$ 
instead of 
$C_0(X^\circ,\End (E|_{X^\circ}))$, $C(X,\End E)$, and $C(\partial X, \End(E|_{\partial X}))$. 
Identifying a continuous function $f$ on $X^\circ$ or $X$ with the 
operator  
$\begin{pmatrix}f&0\\0&0\end{pmatrix}\in \fA$, where $f$ acts by multiplication,
we have natural maps 
$$m_0:C_0(X^\circ)\to \ideal ~\text{ and }~ m:C(X)\to \fA$$ and thus  
a commutative diagram of exact sequences: 
$$
\def\mapup#1{\Big\uparrow\rlap{$\vcenter{\hbox{$\scriptstyle#1$}}$}}
\label{cd}
\begin{array}{ccccccc}
0\longrightarrow&\ideal/\compacts&\longrightarrow&\Boutet/\compacts&{\mathop{\longrightarrow}\limits^{\pi}}&\Boutet/\ideal&\longrightarrow 0
\\                &\mapup{m_0}&&\mapup{m}&&\mapup{b}&
\\0\longrightarrow&C_0(X^\circ)&\longrightarrow&C(X)&{\mathop{\longrightarrow}\limits^{r}}&C(\partial X)&\longrightarrow 0
\end{array} .    
$$
(We do not distinguish between the isomorphic $C^*$-algebras $\Boutet/\ideal$\ and 
the image of $\gamma$.) 
From Lemma \ref{exf} we obtain the commutative grid
 \begin{equation}\label{cde}
    \begin{CD}
      && 0 && 0 && 0\\
      && @AAA @AAA @AAA\\
      0 @>>> C_0(X^\circ) @>>> C(X) @>r>> C(\boundary X) @>>> 0\\
      && @AAA @AAA @AAA\\
      0 @>>> Cm_{0} @>>> Cm @>>> Cb @>>> 0\\
      && @AAA @AAA @AAA\\
      0 @>>> S(\ideal/\compacts) @>>> S(\Boutet/\compacts) @>{S\pi}>>
      S(\Boutet/\ideal) @>>> 0\\
      && @AAA @AAA @AAA\\
      && 0 && 0 && 0
    \end{CD}
  \end{equation}
These six short exact sequences induce six  cyclic long exact sequences in
$K$-theory which we want to analyze next. By Corollary \ref{biso},
$b$ induces an isomorphism in $K$-theory. From Lemma~\ref{l2} and the cyclic
exact sequence of $0\to S(\Boutet/\ideal)\to Cb\to C(\boundary
X)\to 0$ we therefore conclude that $K_*(Cb)=0$. From this in turn we
deduce, using the cyclic exact sequence of $0\to Cm_{0} \to Cm\to
Cb\to 0$, that $Cm_{0}\to Cm$ induces an isomorphism in $K$-theory.

We therefore get the following commutative diagram of cyclic exact
sequences of $K$-theory groups, again using Lemma~\ref{l2} and the natural isomorphisms
$K_{1-i}(SA)\iso K_i(A)$.
\begin{equation}\label{cd6}
  \begin{CD} 
    @>>> K_i(C(X)) @>{m_*}>> K_i(\Boutet/\compacts) @>{\beta}>>K_{1-i}(Cm) @>>>\\
&&@AAA @A{i_*}AA @A{\phi}A{\iso}A\\
    @>>> K_i(C_0(X^\circ)) @>{m_{0*}}>>
    K_i(\ideal/\compacts) @>{\alpha}>> K_{1-i}(Cm_{0}) @>>>
  \end{CD}
\end{equation}

According to Theorem \ref{kergamma}, the principal symbol provides an 
isomorphism $\ideal/\compacts\iso
C_0(S^*X^\circ)$, and $m_0$ becomes the pull back homomorphism $\pi^*$
under this isomorphism.
Since $X$ is connected and $\boundary X\ne \emptyset$, 
there is a nonvanishing section of the cotangent bundle 
(see e.g.~\cite[Proposition 9]{MR1998610} for a proof of 
this well-known fact), which yields a map $C_0(S^*X^\circ)\to
C(X^\circ)$.
Therefore $m_{0*}$ has a split $s$ (which is not necessarily natural). 
Consequently,
$\alpha$ also has a split $s'$. Define now $s'':= i_*\circ s'\circ
\phi^{-1}\colon K_{1-i}(Cm)\to K_i(\Boutet/\compacts)$. An easy
diagram chase shows that $s''$ is a split of $\beta$. Consequently,
our long exact sequence yields  natural short exact sequences
\begin{equation*}
  0\to K_i(C(X)) \xrightarrow{m_*} K_i(\Boutet/\compacts) \xrightarrow{\beta} K_{1-i}(Cm) \to 0,
\end{equation*}
which have a (not necessarily natural) split.
In particular, each element in $ K_i(\Boutet/\compacts)$ can be written as the sum of two elements, 
one in the range of $m_*$ and one in the range of $s''$, thus in the range of $i_*$.

It remains to identify $K_{1-i}(Cm)\iso K_{1-i}(Cm_{0})$. For
this, recall the 
natural 
short exact sequence 
for the ball -- or disc -- completion of the cotangent bundle, extended
to a commutative diagram
\begin{equation}\label{tbs}
  \begin{CD}
    0 @>>>  C_0(T^*X^\circ) @>>> C_0(B^*X^\circ) @>{r}>>
    C_0(S^*X^\circ) @>>> 0\\
    && && @A{\pi^* r_0}A{\sim}A @AA{=}A\\
    && && C_0(B^*X^\circ) @>{\pi^* r_0}>> C_0(S^*X^\circ)\\
    && && @V{r_0}V{\sim}V @VV{=}V\\
    && && C_0(X^\circ) @>{\pi^*=m_{0}} 
                                   >> C_0(S^*X^\circ)
\end{CD}
\end{equation}
Here, $\pi^*$ denotes pull back from the base to the total space of
the bundle, and $r$ and $r_0$ denote restriction to the boundary of the disc
bundle, or the zero section of the disc bundle, respectively; $\sim$ denotes homotopy
equivalences of $C^*$-algebras. Again we have omitted the bundles from the notation.

We get induced short exact mapping cone sequences
\begin{equation}\label{eq:grid}
  \begin{CD}
    0 @>>> SC_0(S^*X^\circ) @>>>  Cr @>>> C_0(B^*X^\circ) @>>> 0\\
    && @AA{=}A @AA{(\pi^*r_0)_*}A @A{\sim}A{\pi^*r_0}A \\
  0 @>>> SC_0(S^*X^\circ) @>>>   C(\pi^* r_0) @>>> C_0(B^*X^\circ)
  @>>> 0\\
  && @VV{=}V  @V{\sim}V{(r_0)_*}V @V{\sim}V{r_0}V\\
 0 @>>> SC_0(S^*X^\circ) @>>> Cm_{0}  
                                    @>>> C_0(X^\circ) @>>> 0
  \end{CD}
\end{equation}

The corresponding cyclic exact $K$-theory sequences together with the
5-lemma imply that the induced maps between the mapping cones induce
isomorphisms in $K$-theory. 

Finally, since $r$ is surjective and $\ker{r}=C_0(T^*X^\circ)$, Lemma~\ref{l3} yields the commutative digram
\begin{equation}\label{cdj}
\def\mapup#1{\Big\uparrow\rlap{$\vcenter{\hbox{$\scriptstyle#1$}}$}}
\begin{array}{ccc}
K_0(SC_0(S^*X^\circ)) &\longrightarrow&K_0(Cr)\\
{\iso\Big\uparrow} &&\iso\mapup{j_*}\\
K_1(C_0(S^*X^\circ))&{\mathop{\longrightarrow}\limits^{\delta}}&K_0(C_0(T^*X^\circ))
\end{array},
\end{equation}
where the lower horizontal arrow is the index mapping for the first row in \eqref{tbs}, and the upper horizontal is induced by 
the first row in \eqref{eq:grid}.

The composition of all these maps gives a natural way to identify
$K_i(Cm)$ with $K_i(C_0(T^*X^\circ))$. This already finishes the proof of Theorem~\ref{t1}. A more detailed explanation of this step will be needed below, in the proof of Theorem~\ref{t2}.

\section{Index Theory}
\label{sec:proof-index-results}

We consider the commutative diagram
\begin{equation}\label{index}
\def\mapdown#1{\Big\downarrow\rlap{$\vcenter{\hbox{$\scriptstyle#1$}}$}}
\def\mapup#1{\Big\uparrow\rlap{$\vcenter{\hbox{$\scriptstyle#1$}}$}}
\begin{array}{ccccccl}
 &K_1(C(X))&{\mathop{\longrightarrow}\limits^{m_*}}           &K_1(\Boutet/\compacts)&\xrightarrow{\beta}&K_0(Cm)         &\\
      &\mapup{} &                                             & \mapup{i_*}           &   &\mapup{\iso}    &\\
 &K_1(C_0(
X^\circ ))&{\mathop{\longrightarrow}\limits^{ {m_{0*}}}}&K_1(\ideal/\compacts) &\xrightarrow{\alpha}&K_0(C m_{0})   & \\
      &                &                              &a\mapup{\iso}          &   &\mapup{\iso}    &\\
       &                &                             &K_0(SC_0(S^* X^\circ ))&\longrightarrow&K_0(C(\pi^*r_0))&\\
        &                &                            &\mapdown{=}             &   &\mapdown{\iso}  &\\
         &                &                           &K_0(SC_0(S^* X^\circ ))&\longrightarrow&K_0(Cr)  &       \\
          &                &                          &c\mapup{\iso}          &   &\mapup{\iso}    &\\
           &                &                         &K_1(C_0(S^* X^\circ ))   &{\mathop{\longrightarrow}\limits^{\delta}}&K_0(C_0(T^* X^\circ ))&\\
            &                &                        &                      &   &\mapdown{\ind_t}&\\
             &                &                       &                      &   &\integers&
\end{array}
\end{equation}                                   
where the first two rows are portions of \eqref{cd6}. The second, third and fourth rows in \eqref{index}\ are portions 
of the cyclic sequences associated to \eqref{eq:grid}\ (notice that, if we use the isomorphism $\ideal/\compacts\iso
C_0(S^* X^\circ )$\ as an identification, then the first column in \eqref{cde}\ is equal to the last row in \eqref{eq:grid}),
while the fourth and fifth rows are just \eqref{cdj}. Note that the
composed
isomorphism $c^{-1}a^{-1}\colon K_1(\ideal/\compacts)\to
K_1(C_0(S^*X^\circ)$ in the left row is exactly the map induced by the
interior symbol.

\begin{definition}
  The map $p$\ in \eqref{eq:main} is the composition of all the maps
  (reverting arrows of isomorphisms when necessary) in the right
  column in \eqref{index}, except $\ind_t$, with the map $\beta$ from
  $K_1(\Boutet/\compacts)$\ to $K_0(Cm)$\ in the first row.
\end{definition}

\begin{remark}
  The definition of $p$ uses the inverse of the
  isomorphism $K_0(Cm_0)\to K_0(Cm)$, which we can not write down
  explicitly -- 
  our argument which proves that the map is an
  isomorphism is actually rather indirect.

  Equivalently, the problem can be restated as replacing a given
  invertible element of 
  $\Boutet/\compacts$ by the sum of elements in the images
  of $m_*$ and $i_*$, respectively, 
  representing the same element in $K_1$.  
  That this is possible is based on the same indirect argument which shows
  that $K_0(Cm_0)\to K_0(Cm)$ is an isomorphism, cf.\ the argument right after 
  \eqref{cd6}. Nevertheless, we will
  see below, in our proof of Fedosov's index formula \eqref{Fedosov},
  that this  representation is actually very useful.
\end{remark}
To prove Theorem \ref{t2}, namely
that
$\ind_t\circ \,p$\ and the Fredholm index are equal on $K_1(\Boutet/\compacts)$, 
it is  enough to show that they are equal on 
the image of $m_*$\ and on the image of $i_*$. 
On the image of $m_*$, both are zero: 
On one hand, the range of $m_*$ consists of equivalence classes 
(modulo  $\compacts$) of invertible multiplication operators. 
Each of these has index zero.  
On the other hand,  the first row in (\ref{index}) is exact, 
thus the range of $m_*$ is mapped to zero. 
The commutativity of \eqref{index} then shows that all we have to
prove is that $\ind_t\circ\,\delta\circ c^{-1}\circ a^{-1}$ is the Fredholm index. 

For that let $\Psi$ denote the $C^*$-closure of the algebra of all
classical pseudodifferential operators of order $0$ on $\widetilde X$ 
in the algebra of all bounded operators on $L^2(\widetilde X)$.
The zero extension on the orthogonal complement of $L^2(X)$ in $L^2(\widetilde X)$ 
defines a *-homomorphism $\ideal_{11}\to\Psi$, where
$\ideal_{11}$ denotes the ideal formed by the   
upper-left corners of $\ideal$. 
That gives us a commutative diagram of exact sequences:
\begin{equation}\label{cdas}
  \begin{CD}
    0 @>>>  \compacts @>>> \ideal_{11} @>>>\ideal_{11}/\compacts @>>> 0\\
    && @VVV @VVV @VV\iota V\\
    0 @>>>  \compacts @>>> \Psi @>>>\Psi/\compacts @>>> 0,
\end{CD}
\end{equation}
where we have denoted by the same symbol the compact ideal in 
the bounded operators on $L^2(X)$ and on $L^2(\widetilde X)$. 
The canonical injection of 
$\ideal_{11}$ into $\ideal$ induces an isomorphism between
$\ideal_{11}/\compacts$ and $\ideal/\compacts$ (see comments right before
Theorem~1 in \cite{MR1998610}). That isomorphism and the naturality of the
index mapping for \eqref{cdas} then imply that
\begin{equation}\label{cdf}
\def\mapdown#1{\Big\downarrow\rlap{$\vcenter{\hbox{$\scriptstyle#1$}}$}}
 \begin{array}{ccc}
K_1(\ideal/\compacts)&\longrightarrow&K_0(\compacts)\iso\integers\\
\mapdown{\iota_*}&&\mapdown{=}\\
K_1(\Psi/\compacts)&\longrightarrow&K_0(\compacts)\iso\integers
\end{array}
\end{equation}
commutes, where the horizontal arrows are the Fredholm-index homomorphisms for $\ideal$ and for $\Psi$.

For any closed manifold, the principal symbol induces an isomorphism
between $\Psi/\compacts$ and the continuous functions on the cosphere 
bundle. This follows from the classical estimate for the norm, 
modulo compacts, of a pseudo-differential operator \cite[Theorem A.4]{KN}. 
We therefore have $K_1(\Psi/\compacts)\iso K_1(C(S^*\widetilde X))$. Modulo this
isomorphism, the Atiyah-Singer index theorem \cite{MR38:5243} states
that the Fredholm-index homomorphism for $\Psi$ is the composition of the topological index 
$\ind_t^{\widetilde X}\colon K_0(C_0(T^*\widetilde X))\to\integers$
with the index mapping for the exact sequence
\[
0\to  C_0(T^*\widetilde X)\to C(B^*\widetilde X) \to C(S^*\widetilde X)\to 0
\]
(see e.g.~\cite[Proposition 15]{MR1998610} for a proof that the classical 
{\em difference bundle} construction indeed gives the
$C^*$-algebra $K$-theory index mapping for this sequence).

Now consider the commutative diagram of exact sequences
\begin{equation}
  \begin{CD}
    0 @>>>  C_0(T^*X^\circ) @>>> C_0(B^*X^\circ) @>>>
    C_0(S^*X^\circ) @>>> 0\\
    && @VVV @VVV @VV\iota V\\
    0 @>>>  C_0(T^*\widetilde X) @>>> C(B^*\widetilde X) @>>>
    C(S^*\widetilde X) @>>> 0.
\end{CD}
\end{equation}
By naturality of the index map, the following diagram commutes:
\begin{equation}\label{fim}
\def\mapdown#1{\Big\downarrow\rlap{$\vcenter{\hbox{$\scriptstyle#1$}}$}}
\begin{array}{ccccc}
K_1(C_0(S^*X^\circ))&{\mathop{\longrightarrow}\limits^{\delta}}&K_0(C_0(T^*X^\circ))\\
\mapdown{\iota_*}&&\mapdown{\tilde\iota}&&\\
K_1(C(S^*\widetilde X))&\longrightarrow&K_0(C_0(T^*\widetilde X))&{\mathop{\longrightarrow}\limits^{\ind_t^{\widetilde X}}}&\integers.
\end{array}
\end{equation}
The Atiyah-Singer index theorem and the commutativity of \eqref{cdf}\ 
imply that the composition of the two lower horizontal and the 
left vertical arrow 
 in (\ref{fim}) gives the Fredholm index, 
hence that
$\ind_t^{\widetilde X}\circ\tilde\iota\circ\delta$
is the Fredholm index for $\ideal$. This proves \eqref{bdmif} since, by definition,
$\ind_t=\ind_t^{\widetilde X}\circ\tilde\iota$. 

To show the last statement, $p={\tt ind}$, look again at
\eqref{index}. As before, it is enough to prove 
that ${\tt ind}\circ i_*=p\circ i_*$ and ${\tt ind}\circ m_*=p\circ
m_*$. Our diagram gives $p\circ i_*=\delta$ and $p\circ m_*=0$,
while ${\tt ind}\circ i_*=\delta$ and ${\tt ind}\circ m_*=0$ are
proven in \cite[Lemmas 16 and 17]{MR1998610}. This shows Theorem
\ref{t2}.

It remains to check the validity of Fedosov's index formula \eqref{Fedosov}.
His proof of homotopy invariance in 
\cite[Proof of Theorem 2.4 in Chapter II]{MR1401125} shows that
the expression only depends on the $K_1$-class represented by the
elliptic operator $A$ in $K_1(\Boutet/\compacts)$, because the formula
is clearly additive for the block sum addition. 
This  
can be considered to be the heart of the proof of the
formula, and we do not offer a new proof for it. 
Once we know that
Fedosov's formula defines a homomorphism from
$K_1(\Boutet/\compacts)$, we can identify it easily with the index
map: It is clear that
the formula is zero for multiplication operators, i.e.~elements in the
image of $K_1(C(X))$. Because of  Theorem \ref{t2} (compare also the
proof above), it suffices to check Formula \eqref{Fedosov} 
for operators $A$ supported in the interior. 
But for those, the formula reduces to the
classical Atiyah-Singer index formula.  
One can use the double of $X$ to get exactly the situation of Atiyah-Singer.

A detail hidden by the simplified notation -- but needed for the 
naturality of the index map -- 
is the existence of a canonical isomorphism
$$K_{0}(C_0(T^*X^\circ, \End(\pi^*E)))\longrightarrow
K_{0}(C_0(T^*X^\circ))$$ 
if $E$ is not the zero bundle.
This is well-known and comes, after a series of standard arguments, from the fact that $E$ is a 
direct summand of a trivial bundle and that the trivial line bundle is
a direct summand in a power of $E$.

\section*{Acknowledgments} Severino Melo was supported by the Brazilian
agency CNPq  (Processos 452780/ 2003-9 and 306214/2003-2). 
Elmar Schrohe had support from the European Research and Training Network
``Geometric Analysis'' (Contract HPRN-CT-1999-0018). 
He thanks Johannes Aastrup and Ryszard Nest for several 
helpful discussions.

{\small

\noindent
Severino T. Melo, Instituto de  Matem{\'a}tica e Estat{\'\i}stica, Universidade de S{\~a}o Paulo, Caixa Postal 66281, 05315-970 S{\~a}o Paulo, email:  {\tt melo@ime.usp.br}\medskip

\noindent
Thomas Schick, Mathematisches Institut, Universit{\"a}t G{\"o}ttingen, Bunsenstr.\ 3-5, \linebreak[4]\mbox{37073 G{\"o}ttingen,}
Germany, email:  {\tt  schick@uni-math.gwdg.de}\medskip

\noindent
Elmar Schrohe, Institut f{\"u}r Mathematik, Universit{\"a}t Hannover, Welfengarten 1,  \linebreak[4] \mbox{30167 Hannover,} Germany, email: {\tt schrohe@math.uni-hannover.de}
}
\end{document}